\documentclass[twoside,12pt,reqno]{amsart}
\usepackage{amsmath,amscd,amssymb,epsfig}
\usepackage{verbatim}

\newcommand{\End}{\operatorname{End}}
\newcommand{\GL}{\operatorname{GL}}
\newcommand{\Id}{\operatorname{Id}}
\newcommand{\str}{\operatorname{str}}
\newcommand{\tr}{\operatorname{tr}}
\newcommand{\sdim}{\operatorname{sdim}}
\newcommand{\qdim}{\operatorname{qdim}}
\newcommand{\sch}{\operatorname{sch}}
\newcommand{\ch}{\operatorname{ch}}
\newcommand{\qd}{\operatorname{\mathsf{d}}}
\newcommand{\e}{\operatorname{e}}
\newcommand{\M}{{\mathcal{M}}}
\newcommand{\T}{{\mathcal{T}}}
\newcommand{\TT}{\mathbb{T}}
\newcommand{\ZZ}{\mathbb{Z}}
\newcommand{\QQ}{\mathbb{Q}}

\newcommand{\CC}{\mathbb{C}}
\newcommand{\NN}{\mathbb{N}}
\newcommand{\sll}{\mathfrak{sl}}
\newcommand{\slto}{\sll(2|1)}
\newcommand{\slmn}{\sll(m|n)}

\newcommand{\osp}{\mathfrak{osp}(2|2n)}
\newcommand{\spn}{\mathfrak{sp}(2n)}
\newcommand{\Uh}{{U_h(\mathfrak{g})}}
\newcommand{\Uhp}{U_h(\mathfrak{g}_{\p0})}

\newcommand{\Uz}{\mathcal{U}_\ZZ}
\newcommand{\qn}[1]{{\left\{#1\right\}}}
\newcommand{\qum}[1]{\widetilde{#1}}
\newcommand{\h}{\ensuremath{\mathfrak{h}}}
\newcommand{\g}{\ensuremath{\mathfrak{g}}}
\newcommand{\borel}{\ensuremath{\mathfrak{b}}}
\newcommand{\nil}{\ensuremath{\mathfrak{n}}}
\newcommand{\p}[1]{\ensuremath{\bar {#1}}}
\newcommand{\D}[1]{\ensuremath{\Delta({#1})} }
\newcommand{\twist}{\ensuremath{\theta}}
\newcommand{\chmap}{\varphi}
\newcommand{\roots}{\Delta}
\newcommand{\weights}{\Lambda}
\newcommand{\chr}{\ZZ[\weights]}
\newcommand{\weyl}{W}
\newcommand{\wta}{\lambda}
\newcommand{\wtaa}{\wta}
\newcommand{\wtb}{\mu}
\newcommand{\wtbb}{\wtb}
\newcommand{\wtc}{\lambda_0}
\newcommand{\wtcc}{\wtc}

\newcommand{\wto}[1]{\lambda_{#1}}
\newcommand{\vz}{\ensuremath{\qum{V}_{\wtcc}}}
\newcommand{\cn}{\xi}
\newcommand{\V}[1]{\widetilde{V}_{#1}}

\newcommand{\bc}{\bar{c}}

\newcounter{bibcount}

\newtheorem{prop}{\bf Proposition}[section]
\newtheorem{defi}[prop]{\bf Definition}
\newtheorem{lem}[prop]{\bf Lemma}
\newtheorem{theo}{\bf Theorem}
\newtheorem{coro}[prop]{\bf Corollary}

\newtheorem{rk}[prop]{Remark}
\newcommand{\epsh}[2]
         {\begin{array}{c} \hspace{-1.3mm}
        \raisebox{-4pt}{\epsfig{figure=#1.eps,height=#2}}
        \hspace{-1.9mm}\end{array}}

\begin{document}
\title[Multivariable link invariants from Lie
  superalgebras]{Multivariable link invariants arising from Lie
  superalgebras of type I}
\author{Nathan Geer}
\address{School of Mathematics\\
Georgia Institute of Technology\\
Atlanta, GA 30332-0160}
\email{geer@math.gatech.edu}
\author{Bertrand Patureau-Mirand}
\address{LMAM, Universit\'e de Bretagne-Sud, BP 573\\
F-56017 Vannes, France }
\email{bertrand.patureau@univ-ubs.fr}
\date{\today}
\begin{abstract} In this paper we construct new links invariants from a type I
  basic Lie superalgebra $\g$.  The construction uses the existence of an
  unexpected replacement of the vanishing quantum dimension of typical module,
  by non-trivial ``fake quantum dimensions.''  Using this, we get a
  multivariable link invariant associated to any one parameter family of
  irreducible $\g$-modules.
\end{abstract}

\maketitle
\setcounter{tocdepth}{1}

\section*{Introduction}

Let $\g$ be a Lie superalgebra of type I, i.e.
$\g$ is equal to $\slmn$ or $\osp$.  Here we assume that $m\neq n$.
Let $r$ be equal to $m+n-1$ if $\g=\slmn$ and $n+1$ if $\g=\osp$.

The quantum dimension associated to a deformed typical $\g$-module is zero.  
This implies that the usual Reshetikhin-Turaev quantum group link invariants
arising from such a module is trivial.  In this paper we show that the usual
Reshetikhin-Turaev construction can be renormalized by non-zero ``fake quantum
dimensions''.  In Section \ref{S:MVA}, we will use this modified construction
to define non-trivial multivariable link invariants.  We will discuss these
multivariable invariants in more detail later in the introduction, now let us 
consider how this modified construction fits into the general theory of
quantum invariants and the representation theory of Lie superalgebras.

The work of this paper has lead to a variety of new mathematical ideas and
relationships.  The authors are working on three subsequent papers which we
will now discuss.  
\begin{itemize}
\item The first of these papers is joint work with V. Turaev.  The paper
  \cite{GPT} will contain a renormalization of the Reshetikhin-Turaev functor
  of a ribbon Ab-category, by ``fake quantum dimensions''.  In the case of
  simple Lie algebras these ``fake quantum dimensions'' are proportional to
  the genuine quantum dimensions.  More interestingly, this paper will contain
  two examples where the genuine quantum dimensions vanish but the ``fake
  quantum dimensions'' are non-zero and lead to non-trivial link invariants.
  The first of these examples recover the hierarchy of invariants defined by
  Akutsu, Deguchi and Ohtsuki \cite{ADO}, using a regularize of the Markov
  trace and nilpotent representations of quantized $\sll(2)$ at a root of
  unity.  These invariants contain Kashaev's quantum dilogarithm invariants of
  knot (see \cite{Kv}).  The second example, is the invariants defined in this
  paper.

  The definition of the ``fake quantum dimensions'' given in \cite{GPT} is
  abstract where the analogous definition in this paper is given by explicit
  formulas.  One can use general theory to show that these definitions are
  equivalent.  The explicit formulas given in this paper are useful when one
  wants to compute the invariant or compare it to other invariants.

\item In the second subsequent paper the authors will use the explicit
  formulas for the ``fake quantum dimensions'' to define ``fake
  superdimensions'' of typical representations of the Lie superalgebra $\g$.
  These ``fake superdimensions'' are non-zero and lead to a kind of supertrace
  on the category representations of $\g$ which is non-trivial and invariant.
  These statements are completely classical statements.  However, the only
  proof we know of uses the quantum algebra and low-dimensional topology
  developed in this paper.

\item We will now discuss the final subsequent paper in relation with the
  multivariable invariants defined in this paper.  In Section \ref{S:MVA} we
  will show that for $c\in \NN^{r-1}$ the pair $(\g,c)$ gives rise to a
  multivariable link invariant $M_{\g}^{c}$.  These invariants associate a
  variable to each component of the link.  There are only a handful of such
  invariants including the multivariable Alexander polynomial and the ones
  defined in \cite{ADO}.  All of these invariants are related to the
  invariants defined in this paper.

  Let us now explain these relationships.  First, in \cite{GP3} we plan on
  showing that the invariant $M_{\sll(m|1)}^{(0,...,0)}$ specializes to the
  multivariable Alexander polynomial.  Second, in order to define their link
  invariants the authors of \cite{ADO} regularize the Markov trace.  Although
  using different methods, the invariants of this paper have a similar
  regularization.  In both case, the standard method using ribbon categories
  or the Markov trace is trivial.  Moreover, both families of invariants are
  generalization of the multivariable Alexander polynomial.  In \cite{GP3}, we
  plan on conjecturing that the invariants $M_{\sll(m|1)}^{(0,...,0)}$, for
  $m\in\NN$, specialize to the hierarchy of invariants defined in \cite{ADO}.
  (Note this specialization depends on $m$ and is different than the
  specialization related to the multivariable Alexander polynomial.)

  Also, in \cite{GP3} we will show that the invariants
  $M_{\sll(m|1)}^{(0,...,0)}$ are related to other invariants.  Let us briefly
  discuss this now.  In \cite{MM}, H. Murakami and J. Murakami show that the
  invariants of \cite{ADO} and the set of colored Jones polynomials have a
  non-trivial intersection.  Moreover, they show that this intersection
  contains Kashaev's invariants.  This result led to a reformulation of
  Kashaev's Volume Conjecture (see \cite{MM}).  In \cite{GP3}, we show that a
  similar result holds for the invariants $M_{\sll(m|1)}^{(0,...,0)}$; namely,
  that the intersection of the set of multivariable link invariants
  $\{M_{\sll(m|1)}^{(0,...,0)}\}_{m\geq 2}$ and the set of colored HOMFLY-PT
  polynomials contains Kashaev's invariants.  We also show that the invariants
  $M_{\sll(m|1)}^{(0,...,0)}$ are multivariable generalization of the set of
  two variable invariants defined by Links and Gould \cite{Dw,LG}.
\end{itemize}

Lie superalgebras have previously been used to construct invariants with more
than one variable.  For example see \cite{GP,LG,LGZ,LZ}.  It should be noted
that in \cite{LG,LGZ} the use of cutting one strand of a link in order to
avoid the vanishing of the quantum dimension has already been used.  In these
papers all the strands of a braid are colored with the same module and the
authors use the Markov trace to construct a non-trivial invariant.  Here we
use links whose components are colored with different modules and work with
ribbon functors.  As explained above, this approach gives the potential for
new constructions in low-dimensional topology and applications in
representation theory.

We now state the main results of this paper more precisely.  First, we recall
some results from the theory of Lie superalgebras (for more details see
Section \ref{s:prel}).  Every irreducible finite-dimensional $\g$-module has a
highest weight $\lambda\in\h^*$ (where $\h$ is the Cartan sub-superalgebra).
Moreover, the set of isomorphism classes of irreducible finite-dimensional
$\g$-modules are in one to one correspondence with the set of dominant
weights.  These modules are parameterized by $\NN^{r-1}\times\CC$ and are
divided into two classes: typical and atypical.  Each highest weight
$\g$-module $V$ can be deformed to a highest weight topologically free
$\Uh$-module $\qum{V}$, where $\Uh$ is the Drinfeld-Jimbo superalgebra
associated to $\g$.  We say $\qum{V}$ is a typical $\Uh$-module if $V$ is a
typical $\g$-module.

Let $F$ be the usual Reshetikhin-Turaev functor from the category of framed
tangles colored by topologically free $\Uh$-modules of finite rank to the
category of $\Uh$-modules (see \cite{Tu}).  In Section \ref{s:proofF'}, we
define a map $\qd$ from the set of typical $\Uh$-module to the ring
$\CC[[h]][h^{-1}]$.  If $T_{\lambda}$ is a framed $(1,1)$-tangle colored by
$\Uh$-modules such that the open string is colored by the deformed typical
module $\qum{V}(\lambda)$ of highest weight $\lambda$, then
$F(T_{\lambda})=x.\Id_{\qum{V}(\lambda)}$, for some $x$ in $\CC[[h]]$.  We set
$F'(T_{\lambda})=x.\qd(\lambda)$.

If $L$ is the colored link given by the closure of $T_{\lambda}$, we will see
at the end of Section \ref{s:prel} that $F(L)=x.\qdim(\qum{V}(\lambda))=0$.
For this reason we think of $\qd$ as a replacement for the quantum dimension
$\qdim$.  In Section \ref{s:proofF'}, we show that this regularization makes
the map $F'$ into a well defined framed colored link invariant.  In
particular, we prove the following theorem.
\begin{theo} \label{thF'} The map $F'$ induces a well defined invariant of
  framed links colored by at least one typical $\Uh$-module.  In other words,
  if $L$ is a framed link colored by $\Uh$-modules and the closure of
  $T_{\lambda}$ is equal to $L$ then the map given by $L\mapsto
  F'(T_{\lambda})$ is a well defined framed colored link invariant.
\end{theo}
In Section \ref{S:MVA} we will use $F'$ to show that there exists
multivariable link invariants $M^{c}_{\slmn}$ and $M^{c}_{\osp}$ for each
$m,n\in\NN^*$ and $c\in \NN^{r-1}$.  Let us now make this statement precise.

Let $c\in \NN^{r-1}$.  Then for all but a finite number of $a\in\CC$ the
module corresponding to $(c,a)$ is typical.  Let $\TT_c$ be this set of
complex numbers.  Let $\wto a^c \in\h^{*}$ be the weight corresponding to
$(c,a)$.  The proof of the following theorem can be found in Section
\ref{S:MVA}.

\begin{theo}\label{T:Mpoly2}
  Let $L'$ be a framed link with $k$ ordered components. Let $L$ be the
  non-framed link which underlies $L'$.  For each $c\in \NN^{r-1}$, there
  exists an multivariable link invariant $M_{\g}^{c}$ with the following
  properties.
  \begin{enumerate}
  \item \label{I:Mpoly1} If $k=1$ then $M_{\g}^{c}(L)$ takes values in
    $(M^c_1(q,q_1))^{-1} \ZZ[q^{\pm1},q_1^{\pm1}]$,
  \item If $k\geq 2$ then $M_{\g}^{c}(L)$ takes values in
    $\ZZ[q^{\pm1},q_1^{\pm1},\ldots,q_k^{\pm1}]$,
  \item If $(\cn_1,\ldots,\cn_k)\in(\TT_c)^k$ and the $i$th component of $L'$
    is colored by the typical module of weight $\wto{\cn_i}^{c}$ then
    $$F'(L')=\e^{\sum lk_{ij}<\wto{\cn_i}^{c},\wto{\cn_j}^{c}+2\rho>h/2}
    {M_{\g}^{c}(L)|_{q_i=\e^{\cn_ih/2}}}$$
  \end{enumerate}
  where $(lk_{ij})$ is the linking matrix of $L'$; the bilinear form $<.,.>$,
  and element $\rho\in \h $ are defined in subsection \ref{SS:Liesuper}; and
  $M^c_1$ is defined in Lemma \ref{L:Mc}.
\end{theo}
More generally we define in Theorem \ref{T:MCpoly2} (page \pageref{T:MCpoly2})
a $(k+1)$-multi\-variable invariant of ordered links whose $k$ components are
colored by different elements of $\NN^{r-1}$.

In \cite{GP} the authors give proves of Theorem \ref{thF'} and Theorem
\ref{T:Mpoly2} in the case of $\slto$ and $c=(0,\ldots,0)$.  Many of the
results, in \cite{GP}, are proved by calculations made by hand.  In this paper
our proofs are based on more general techniques rooted in the representation
theory of $\g$.

Still these techniques require us to assume that $\g$ is a Lie superalgebras
of type I.  We make this assumption because the character formulas for typical
modules of other basic classical Lie superalgebras (i.e. not of type I), are
more complicated.  Also, basic classical Lie superalgebras of type II, do not
have one-parameter families of modules and thus no natural way to construct
multivariable invariants.  However, for these Lie superalgebras one should
still be able to define a map $\qd$ and thus a framed colored link invariant.
\subsection*{Acknowledgments}  N.G. would like to thank LMAM,
Universit\'e de Bretagne-Sud for their generous hospitality.
\section{Preliminaries}\label{s:prel}

In the section we review background material that will be used in the
following sections.

A \emph{super-space} is a $\ZZ_{2}$-graded vector space $V=V_{\p 0}\oplus
V_{\p 1}$ over $\CC$.  We denote the parity of a homogeneous element $x\in V$
by $\p x\in \ZZ_{2}$.  We say $x$ is even (odd) if $x\in V_{\p 0}$ (resp.
$x\in V_{\p 1}$).  A \emph{Lie superalgebra} is a super-space $\g=\g_{\p 0}
\oplus \g_{\p 1}$ with a super-bracket $[\: , ] :\g^{\otimes 2} \rightarrow
\g$ that preserves the $\ZZ_{2}$-grading, is super-antisymmetric
($[x,y]=-(-1)^{\p x \p y}[y,x]$), and satisfies the super-Jacobi identity (see
\cite{K}).  Throughout, all modules will be $\ZZ_{2}$-graded modules (module
structures which preserve the $\ZZ_{2}$-grading, see \cite{K}).

\subsection{Lie superalgebras of type I}\label{SS:Liesuper}
In this subsection we recall notation and properties related to Lie
superalgebras of type I and modules over such Lie superalgebras.  Modules over
Lie superalgebras of type I are different in nature than modules over
semi-simple Lie algebras.  For example, each Lie superalgebra of type I has
one parameter families of modules.

Let $\g=\g_{\p 0}\oplus \g_{\p 1}$ be a Lie superalgebra of type I, i.e.
$\g$ is equal to $\slmn$ or $\osp$.  We will assume that $m\neq n$.  Let
$\borel$ be the distinguished Borel sub-superalgebra of $\g$.  Then $\borel$
can be written as the direct sum of a Cartan sub-superalgebra $\h$ and a
positive nilpotent sub-superalgebra $\nil_{+}$.  Moreover, $\g$ admits a
decomposition $\g=\nil_{-}\oplus \h \oplus \nil_{+}$.  Let $\weyl$ be the Weyl
group of the even part $\g_{\p0}$ of $\g$.

Let $\roots_{\p 0}^{+}$ (resp. $\roots_{\p 1}^{+}$) be the even (resp. odd)
positive roots.  Let $\rho_{\p0}$ (resp. $\rho_{\p1}$) denote the half sum of
all the even (resp. odd) positive roots.  Set $\rho= \rho_{\p0}-\rho_{\p1}$.
A positive root is called \emph{simple} if it cannot be decomposed into a sum
of two positive roots.

A Cartan matrix associated to a Lie superalgebra is a pair consisting of a $r
\times r$ matrix $A=(a_{ij})$ and a set $\tau\subset \{1,\ldots,r\}$
determining the parity of the generators.  Let $(A,\tau)$ be the Cartan
arising from $\g$ and the distinguished Borel sub-superalgebra $\borel$.  Here
the set $\tau=\{s\}$ consists of only one element because of our choice of
Borel sub-algebra $\borel$. (See the appendix.)

There are $d_{1},\ldots,d_{r}$ in $\{\pm1,\pm2\}$ such that the matrix
$(d_{i}a_{ij})$ is symmetric. Here we will assume $d_{1}=1$.  Let $<.,.>$ be
the symmetric non-degenerate form on $\h$ determined by
$<h_{i},h_{j}>=d_{j}^{-1}a_{ij}$.  This form gives an identification of $\h$
and $\h^{*}$.  Moreover, the form $<.,.>$ induces a $\weyl$-invariant bilinear
form on $\h^{*}$, which we will also denote by $<.,.>$.

By Proposition 1.5 of \cite{Kac78} there exists $e_{i}\in \nil_{+}$, $f_{i}\in
\nil_{-}$ and $h_{i}\in\h$ for $i=1,\ldots,r$ such that the Lie superalgebra
$\g$ is generated by $e_{i}, f_{i}, h_{i}$ where
\begin{align*}
[e_{i},f_{j}]=& \delta_{ij}h_{i}, & [h_{i},h_{j}]=&0, &
[h_{i},e_{j}]=&a_{ij}e_{j}, & [h_{i},f_{j}]=&-a_{ij}f_{j}.
\end{align*}

Let $\wta\in\h^{*}$ be a linear functional on $\h$.  Kac \cite{K} defined a
$\g$ irreducible highest weight module $V(\wta)$ of weight $\wta$ with a
highest weight vector $v_{0}$ having the property that $h.v_{0}=\wta(h)v_{0}$
for all $h\in\h$ and $\nil_{+}v_{0}=0$.  Let $a_{i}=\wta(h_{i})$.  In \cite{K}
Kac showed that $V(\wta)$ is finite-dimensional if and only if $a_{i}\in \NN$
for $i\neq s$.  Therefore, $a_{s}$ can be an arbitrary complex number.
Irreducible finite-dimensional $\g$-modules are divided into two classes
typical and atypical.

There are many equivalent definitions for a weight module to be typical (see
\cite{Kac78}).  Here we say that $V(\wta)$ is typical if it splits in any
finite-dimensional $\g$-module.  By Theorem 1 of \cite{Kac78} this is
equivalent to requiring that
\begin{equation}
  \label{E:typ}
  <\wta+\rho,\alpha>\neq 0
\end{equation}
for all $\alpha \in \roots_{\p 1}^{+}$.
If $V(\wta)$ is typical we will say the weight $\wta$ is typical.  Also note
that typical weights are dense in the space of weights corresponding to
finite-dimensional modules.  In particular, if $a_{i}\in \NN$ for $1\leq i\leq
r$ and $i\neq s$ then there are only finitely many atypical weights with
$a_{i}=\wta(h_{i})$.  Furthermore, if $\wta$ is atypical then
$a_s=\wta(h_{s})\in\ZZ$.

Let $\weights\simeq \ZZ^{r-1}\times\CC$ be the group (or weight ``lattice'')
of weights taking integer values on $h_i$ for $i\neq s$.  For $\alpha\in
\weights$ we denote its image in $\chr$ by $e^{\alpha}$.

Let $L'_0$, $L'_1$ and $L_1$ be the following elements of the ring $\chr$ of
characters:
$$L'_{0}=\prod_{\alpha \in
  \roots_{\p0}^{+}}(\e^{\alpha/2}-\e^{-\alpha/2})\quad,\quad L_1=\prod_{\alpha
  \in \roots_{\p1}^{+}} (\e^{\alpha/2}+\e^{-\alpha/2})$$
$$\quad\text{ and }L'_{1}=\prod_{\alpha \in
  \roots_{\p1}^{+}}(\e^{\alpha/2}-\e^{-\alpha/2}).$$
If $V(\wta)$ is typical then from $(2.2')$ of \cite{Kac78} we have the
following formula for the (super-) character of $V(\wta)$:
\begin{equation}
  \label{E:sch}
  \begin{split}
    \sch(V(\wta))=(L'_1/L'_0)\sum_{w\in W}\epsilon(w)\e^{w(\wta+\rho)}\\
    \ch(V(\wta))=(L_1/L'_0)\sum_{w\in W}\epsilon(w)\e^{w(\wta+\rho)}
  \end{split}
\end{equation}
with $\epsilon: \weyl \rightarrow \{\pm 1\}$ the function given by
$\epsilon(w)=(-1)^{c}$ where $c$ is the number of reflections in the
expression of $w$.

\begin{prop}\label{P:formSch}
  The super-character (resp. the character) of a typical $\g$-module $V(\wta)$
  has the form $\sch(V(\wta))=\chi'_{1}\chi_{0}(\wta)$ (resp.
  $\ch(V(\wta))=\chi_{1}\chi_{0}(\wta)$) where
  $$\chi'_1=\prod_{\alpha \in
    \roots^{+}_{\p1}}(1-\e^{-\alpha})\quad,\quad \chi_1=\prod_{\alpha \in
    \roots^{+}_{\p1}}(1+\e^{-\alpha})$$
  and $\chi_{0}(\wta)$ is the character of the even finite-dimensional
  irreducible $\g_{\p0}$-module with highest weight $\wta$.
\end{prop}
\begin{proof}
  Let $V(\wta)$ be a typical $\g$-module.  From \cite{Kac78} Proposition
  1.7(c) we have that $w(\rho_{1})=\rho_{1}$ for all $w\in W$ and so equation
  \eqref{E:sch} can be rewritten as
  \begin{align}
    \sch(V(\wta))&=(L'_1/L'_0)\sum_{w\in
      W}\epsilon(w)\e^{w(\wta+\rho_{0})-\rho_{1}} \notag \\
    &=\prod_{\alpha \in \roots_{\p1}^{+}}(1-\e^{-\alpha})\frac{\sum_{w\in
        W}\epsilon(w)\e^{w(\wta+\rho_{0})}}{\prod_{\alpha \in
        \roots_{\p0}^{+}}(\e^{\alpha/2}-\e^{-\alpha/2})}\label{E:sch2}
  \end{align}
  where the last equality follows from the fact that $\rho_{1}$ is the half
  sum of the positive odd roots.  Now the fraction in equation \eqref{E:sch2}
  is the character of the even finite-dimensional irreducible
  $\g_{\p0}$-module with highest weight $\wta$.  A similar argument shows that
  $\ch(V(\wta))=\chi_{1}\chi_{0}(\wta)$.
\end{proof}

\begin{rk}
  In formulas \eqref{E:typ}-\eqref{E:sch2} we use the fact that $\g$ is a Lie
  superalgebras of type I.  In particular, in the notation of \cite{Kac78} we
  have $\roots_{\p0}^{+}=\overline{\roots_{\p0}^{+}}$ and
  $\epsilon'=\epsilon$.
\end{rk}

\begin{rk}
  Here we only consider ``even'' irreducible modules: $V(\wta)$, i.e. modules
  with an even highest weight vector.  Every such irreducible module has an
  ``odd'' analog $V(\wta)^-$ obtained by just taking the opposite parity.
  $V(\wta)^-$ is isomorphic with $V(\wta)$ with an odd isomorphism.  Remark
  that the tensor product of two ``even'' modules may contain ``odd'' modules.
  The character and super-character of an ``odd'' module is
$$\ch(V(\wta)^-)=\ch(V(\wta))\quad\text{and}\quad
\sch(V(\wta)^-)=-\sch(V(\wta)).$$
\end{rk}


\subsection{The quantization $\Uh$}\label{SS:quant}

Let $h$ be an indeterminate.  Set $q=\e^{h/2}$.  We adopt the following
notations:
$$q^z=\e^{zh/2}\quad\textrm{and}\quad\qn z=q^z-q^{-z}.$$
\begin{defi}[\cite{Yam94}]\label{D:Usl}
  Let $\g$ be a Lie superalgebra of type I.  Let $(A,\{s\})$ be the Cartan
  matrix arising from the distinguished Borel sub-super\-algebra (see section
  \ref{SS:Liesuper}).  Let $\Uh$ be the $\CC[[h]]$-Hopf superalgebra generated
  by the elements $h_{i},E_{i}$ and $ F_{i}, $ $i = 1\cdots r$, satisfying the
  relations:
  \begin{align*}
    [h_{i},h_{j}] &=0, & [h_{i},E_{j}]=&a_{ij}E_{j}, &
    [h_{i},F_{j}]=&-a_{ij}F_{j},
  \end{align*}
  \begin{align*}
    [E_{i},F_{j}]=&\delta_{i,j}\frac{q^{h_{i}}-q^{-h_{i}}}{q-q^{-1}}, &
    E_{s}^{2}=&F_{s}^{2}=0,
  \end{align*}
  plus the quantum Serre-type relations (see Definition 4.2.1 of
  \cite{Yam94}).  Here $[,]$ is the super-commutator given by
  $[x,y]=xy-(-1)^{\p x \p y}yx$.  All generators are even except for $E_{s}$
  and $F_{s}$ which are odd.  The coproduct, counit and antipode are given by
  \begin{align*}
    \label{}
    \D {E_{i}}= & E_{i}\otimes 1+ q^{-h_{i}} \otimes E_{i}, &
    \epsilon(E_{i})= & 0 & S(E_{i})=&-q^{h_{i}}E_{i}\\
    \D {F_{i}} = & F_{i}\otimes q^{h_{i}}+ 1 \otimes F_{i}, &
    \epsilon(F_{i})= &0 & S(F_{i})=&-F_{i} q^{-h_{i}}\\
    \D {h_{i}} = & h_{i} \otimes 1 + 1\otimes h_{i}, & \epsilon(h_{i}) = & 0 &
    S(h_{i})= &-h_{i}.
  \end{align*}
\end{defi}

Khoroshkin, Tolstoy \cite{KTol} and Yamane \cite{Yam94} showed that $\Uh$ has
an explicit $R$-Matrix $R$.  We will now recall their results.  Let
$\exp_{q}(x):= \sum_{n=0}^{\infty}x^{n}/(n)_{q}!$ be the ``q-exponential,''
where $(n)_{q}!:=(1)_{q}(2)_{q}\cdots (n)_{q}$ and $(k)_{q}:=(1-q^{k})/(1-q)$.
Fix a normal ordering on the set of positive roots $\roots^{+}$.  Let
$E_{\alpha}$ and $F_{\alpha}$ for $\alpha \in \roots^{+}$ be the $q$-analogs
of the Cartan-Weyl generators where $E_{\alpha_{i}}=E_{i}$ and
$F_{\alpha_{i}}=F_{i}$ for any simple root $\alpha_{i}$ (see Section 3 of
\cite{KTol}).  Let $\{e_{\alpha},f_{\alpha},h_{i}:\alpha\in
\roots^{+},i\in\{1,...,r\}\}$ be the Cartan-Weyl basis of $\g$.  For
$\alpha\in \roots^{+}$ let $a_{\alpha}$ be the function of $q$ defined by
$$[E_{\alpha},
  F_{\alpha}]=a_{\alpha}(q^{h_{\alpha}}-q^{-h_{\alpha}})/(q-q^{-1})$$
  where $h_{\alpha}:=[e_{\alpha},f_{\alpha}]$.
Let
$$\check{R}_{\alpha}:=\exp_{q}\big((-1)^{\p
  \alpha}a_{\alpha}^{-1}(q-q^{-1})(E_{\alpha}\otimes F_{\alpha})\big).$$ 
Also, let
\begin{equation}
  \label{E:Rcheck}
  \check{R}=\prod_{\alpha \in\roots^{+}}\check{R}_{\alpha},
\end{equation}
where the order in the product is given by the chosen fixed normal ordering.
Let $(d_{ij})$ be the inverse of the matrix $(a_{ij}/d_j)$.  Set
\begin{equation}
  \label{E:FormOfK}
  K=q^{\sum_{i,j}^{r}d_{ij}h_{i}\otimes h_{j}}.
\end{equation}
Then the $R$-matrix is of the form $R=\check{R}K$.  Remark for $\alpha\in
\roots^{+}_{\p1}$ it follows from Zhang \cite{ZA,ZC} that $a_{\alpha}=1$.  One
can also see from Yamane \cite{Yam94} that the $R$-matrix has the above form.

We say a $\Uh$-module $W$ is topologically free of finite rank if it is
isomorphic as a $\CC[[h]]$-module to $V[[h]]$, where $V$ is a
finite-dimensional $\g$-module.  Let $\M$ be the category of topologically
free of finite rank $\Uh$-modules.  A standard argument shows that $\M$ is a
ribbon category (for details see \cite{G04B}).  Let $V,W $ be objects of $\M$.
We denote the braiding and twist morphisms of $\M$ as
\begin{align*}
  c_{V,W}:&V\otimes W \rightarrow W \otimes V, & \twist_{V}:&V\rightarrow V
\end{align*}
respectively.  We also denote the duality morphisms of $\M$ as
\begin{align*}
\label{}
b_{V} :&\CC[[h]]\rightarrow V\otimes V^{*}, & d_{V}': & V\otimes
V^{*}\rightarrow \CC[[h]]
\end{align*}

Let $\T=Rib_\M$ be the ribbon category of ribbon graphs over $\M$ in the sense
of Turaev (see \cite{Tu} chapter I) where coupons are colored by even
morphisms.  The set of morphisms $\T((V_1,\ldots,V_n),(W_1,\ldots,W_m))$ is a
space of formal linear combinations of ribbon graphs colored by objects of
$\M$.  Let $F$ be the usual ribbon functor from $\T$ to $\M$ (see \cite{Tu}).

In \cite{G05}, Geer shows that a weight $\g$-module $V(\wta)$ can be deformed
to a weight $\Uh$-module $\qum{V}(\wta)$.  In particular, Geer shows that the
characters of $V(\wta)$ and $\qum{V}(\wta)$ are equal and that as a
super-space $\qum{V}(\wta)$ is equal to $V(\wta)[[h]]$.

We say that $V\in\M$ is \textit{irreducible} if $\End_\Uh(V)=\CC[[h]].\Id_V$.
Then the deformation $\qum{V}(\wta)$ is irreducible for every
finite-dimensional irreducible weight $\g$-module $V(\wta)$.

It is well known that the super-dimension of any typical $\g$-module is zero.
Then an argument using the Kontsevich integral shows that the quantum
dimension of any deformed typical $\Uh$-module factors as $x.\sdim(V)$ (for
some $x\in \CC[[h]]$) and thus is zero.  It follows that the functor $F$ is
zero on all closed ribbon graph over $\M$ with at least one color which is a
deformed typical module.  For this reason it can be difficult to construct
nontrivial link invariants from typical $\g$-modules.
\section{Proof of Theorem \ref{thF'}}\label{s:proofF'}

In this section we define the map $\qd$ and prove a series of lemmas which
lead to the proof of Theorem~\ref{thF'}.
\begin{defi}
  If $T\in\T\left(V,V\right)$ where $V\in\M$ is irreducible then
  $F(T)=x.\Id_{V} \in \End_{\Uh}(V)$ for some $x\in\CC[[h]]$.  We define the
  bracket of $T$ to be
  $$<T>=x.$$
\end{defi}
For example, if $V, V'$ are modules of $\M$ such that $V'$ is irreducible, we
define
$$S'(V,V')=\left< \epsh{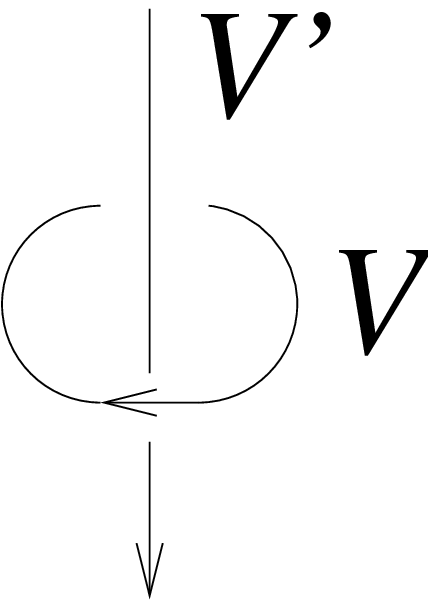}{8ex}\right>.$$
When $V=\qum{V}(\wta)$ and $V'=\qum{V}(\wtb)$ are irreducible highest weight
modules with highest weights $\wta$ and $\wtb$ we write $S'(\wta,\wtb)$ for
$S'(V,V')$.

For any weight $\beta$, let $\chmap_{\beta}: \chr \rightarrow \CC[[h]]$ be
given by
\begin{equation}
  \label{E:chmap}
  \e^{\alpha} \mapsto q^{2<\alpha, \beta>}.
\end{equation}
\begin{prop} \label{P:S'}
  $$S'(\wta,\wtb)=\chmap_{\wtb+\rho}(\sch(V(\wta)).$$
\end{prop}
\begin{proof}

  Let $(v_{i})$ be a basis of $V(\wta)$ such that $v_{i}$ is a weight vector
  of weight $\alpha_{i}\in \h^{*}$.

  Let $w_{\wtb}$ be a highest weight vector of $\qum{V}(\wtb)$.  Recall the
  $R$-matrix is of the form $R=\check{R}K$.  We will use equation
  \eqref{E:FormOfK} to show
  \begin{align}
    \label{E:K1}
    K(w_{\wtb}\otimes v_{i}) &=q^{<\wtb,\alpha_{i}>}w_{\wtb}\otimes v_{i}&
    K(v_{i}\otimes w_{\wtb}) &=q^{<\wtb,\alpha_{i}>}v_{i}\otimes w_{\wtb}.
 \end{align}
 Indeed, let $v$ and $v'$ be two weight vectors of respective weight $\nu$ and
 $\nu'$.  As $(a_{ij}/d_j)$ is the matrix of $<.,.>_{|\h\otimes\h}$ in the
 basis $(h_i)$ and so its inverse $(d_{ij})$ is the matrix of
 $<.,.>_{|\h^*\otimes\h^*}$ in the dual basis.  Now the scalars
 $x^{(_{'})}_i=\nu^{(_{'})}(h_i)$ by which $h_i$ acts on $v^{(_{'})}$ are
 precisely the coordinates of $\nu^{(_{'})}$ in the basis dual to $(h_i)$.
 Hence $\sum d_{ij}h_i\otimes h_j (v\otimes v')=(\sum
 x_id_{ij}x'_j)v\otimes v'=<\nu,\nu'>v\otimes v'$.\\

 We now give two facts.  Let $v$ be any weight vector of $\qum{V}(\wta)$ of
 weight~$\eta$.
 \begin{description}
 \item[Fact 1] $R(w_{\wtb}\otimes v)=q^{<\wtb,\eta>}(w_{\wtb}\otimes
   v)$.\\
   This fact follows from equations \eqref{E:Rcheck} and \eqref{E:K1} and the
   property that $E_{\alpha}w_{\wtb}=0$ for $\alpha\in \roots^+$.
 \item[Fact 2] All the pure tensors of the element $(\check{R}-1)(v\otimes
   v_{\wtb})\in \qum{V}(\wta)\otimes \qum{V}(\wtb)$ are of the form $v'
   \otimes w'$ where $w'$ is a weight vector of $\qum{V}(\wtb)$ and $v'$ is a
   weight vector of $\qum{V}(\wta)$ whose weight is of strictly higher
   order than that of the weight of $v$.\\
   Fact 2 is true because $E_{\alpha}^{n}v$ (for $\alpha\in \roots^+$ and
   $n\in \NN^*$) is zero or a weight vector whose weight is of strictly higher
   order than the weight of $v$.
 \end{description}

 We will now compute $S'(\wta,\wtb)$ directly.  Let $V$ be an object of $\M$
 and recall that the duality morphisms $b_{V} :\CC[[h]]\rightarrow V\otimes
 V^{*}$ and $d_{V}': V\otimes V^{*}\rightarrow \CC[[h]]$ are defined as
 follows.  The morphism $b_{V}$ is the $\CC[[h]]$-linear extension of the
 coevaluation map on the underlying $\g$-module.  In particular,
 \begin{equation}
   \label{E:bv}
   b_{\qum{V}(\wta)}(1)=\sum v_{i}\otimes v_{i}^{*}
 \end{equation}

 As in the case of semi-simple Lie algebras we have
 \begin{equation}
   \label{E:K2}
   d_{\qum{V}(\wta)}' (v\otimes f)= (-1)^{\p v \p
     f}f(q^{2<\eta,\rho>}v)
 \end{equation} where $v$ is a weight vector of
 $\qum{V}(\wta)$ of weight $\eta\in\h^{*}$.  Note in the above equation
 $q=\e^{h/2}$.  Consider the element $S\in \End_{\Uh}(\qum{V}(\wtb))$
 given by
 $$(\Id_{\qum{V}(\wtb)}\otimes
 d_{\qum{V}(\wta)}')\circ(c_{\qum{V}(\wta),\qum{V}(\wtb)}\otimes
 \Id_{\qum{V}(\wta)^{*}})\circ ( c_{\qum{V}(\wtb),\qum{V}(\wta)}\otimes
 \Id_{\qum{V}(\wta)^{*}})\circ(\Id_{\qum{V}(\wtb)}\otimes
 b_{\qum{V}(\wta)}).$$ 
 To simplify notation set $S=(X_{1})(X_{2})(X_{3})(X_{4})$ where $X_{i}$ is
 the corresponding morphism in the above formula.  The morphism $S$ is
 determined by its value on the highest weight $w_{\wtb}$.  By definition
 $S(w_{\wtb})=S'(\wta,\wtb)w_{\wtb}$, so it suffices to compute $S(w_{\wtb})$.
 \begin{align}
   \label{E:valueS}
   S(w_{\wtb})&=(X_{1})(X_{2})(X_{3})\Big(w_{\wtb}\otimes
   \sum_{i}(v_{i}\otimes v_{i}^{*})\Big)\notag\\
   &= (X_{1})(X_{2})\Big( \sum_{i} q^{<\wtb,\alpha_{i}>}
   v_{i}\otimes w_{\wtb}\otimes v_{i}^{*}\Big)\notag \\
   &= (X_{1})\Big( \sum_{i} \big( q^{2<\wtb,\alpha_{i}>} w_{\wtb}\otimes
   v_{i}\otimes v_{i}^{*}\big) + \sum_{k} w_{k}' \otimes
   v_{k}' \otimes z_{k} \Big) \notag\\
   &=\sum_{i}(-1)^{\p {v_i}} q^{2<\wtb +\rho, \alpha_{i}>}w_{\wtb}
 \end{align}
 where $z_{k}= v_{i}^{*}$ (for some $i$), $v'_{k}$ is a weight vector of
 $\qum{V}(\wta)$ whose weight is of strictly higher weight than that of the
 weight of $z_{k}^{*}$ and $w'_{k}$ is a weight vector of $\qum{V}(\wtb)$.
 Moreover, the first equality of the above equation follows from \eqref{E:bv},
 the second from Fact 1, the third from \eqref{E:K1} and Fact 2, and finally
 the fourth from \eqref{E:K2} and the fact that $z_{k}(v'_{k})=0$.  The key
 observation in this proof is that Facts 1 and 2 imply that in the above
 computation the only contribution of the action of the $R$-matrix comes from
 $K$.

 Since the super-character of $V(\wta)$ is equal to $\sum
 (-1)^{\p\alpha_{i}}\e^{\alpha_{i}}$ the proposition follows from equation
 \eqref{E:valueS}.
\end{proof}
%
\begin{coro}\label{C:SnotZero}
  Let $\wta$ be a typical weight. Then $S'(\wta,\wtb)=0$ if and only if $\wtb$
  is atypical.
\end{coro}
\begin{proof}
  From Proposition \ref{P:formSch} we have
  $$S'(\wta,\wtb)= \prod_{\alpha \in
    \roots^{+}_{\p1}}\chmap_{\wtb+\rho}(1-\e^{-\alpha})
  \chmap_{\wtb+\rho}(\chi_{0}(\wta)).$$ 
  Since $\chi_{0}(\wta)$ is a character of an even finite-dimensional
  $\g_{\p0}$-module we have $\chmap_{\wtb+\rho}(\chi_{0}(\wta))\neq 0$ for any
  weight $\wtb$.  Now by the definition of a typical module (see equation
  \eqref{E:typ}) we have $ \prod_{\alpha \in
    \roots^{+}_{\p1}}\chmap_{\wtb+\rho}(1-\e^{-\alpha})$ is non-zero if and
  only if $\wtb$ is typical.
\end{proof}

Recall the definitions of $\chmap_{\wta}$ and $L'_{i}$ given in equation
\eqref{E:chmap} and Subsection \ref{SS:Liesuper}, respectively.
\begin{lem}\label{L:def}
  Let $\wta$ be a typical weight.  Then $\frac{\chmap_{\wta+\rho}(L'_{0})}
  {\chmap_{\wta+\rho}(L'_{1})\chmap_{\rho}(L'_{0})}$ is an element of
  $h^{-k_{\g}}\CC[[h]]$ where $k_{\g}$ is the number of odd positive roots,
  i.e. $k_{\g}=|\roots_{\p1}^{+}|$.
\end{lem}
\begin{proof}
  We have $L'_{1}=\e^{\rho_{1}}\prod_{\alpha \in
    \roots^{+}_{\p1}}(1-\e^{-\alpha})$.  Therefore, since $\wta$ is typical
  then by equation \eqref{E:typ} we have $\chmap_{\wta+\rho}(L'_{1})$ is
  non-zero element of $\CC[[h]]$.  Moreover, since the product in $L'_{1}$ is
  taken over $\roots_{\p1}^{+}$ we have
  $h^{-|\roots_{\p1}^{+}|}\chmap_{\wta+\rho}(L'_{1})$ is invertible in
  $\CC[[h]]$.  The lemma follows from the fact that
  $\chmap_{\wta+\rho}(L'_{0})/\chmap_{\rho}(L'_{0})\in \CC[[h]]$.
\end{proof}
\begin{defi}\label{D:d}
  Let $\wta$ be a typical weight.  Define
  $$\qd(\qum V(\wta))=\qd(\wta)=
  \frac{\chmap_{\wta+\rho}(L'_{0})}
  {\chmap_{\wta+\rho}(L'_{1})\chmap_{\rho}(L'_{0})} \in h^{-k_{\g}}\CC[[h]]$$
  where we use the notation of Lemma \ref{L:def}.
\end{defi}

\begin{lem}\label{L:sym}
  Let $\wta$ and $\wtb$ be typical weights.  Then
  $$\qd(\wtb)S'(\wta,\wtb)=\qd(\wta)S'(\wtb,\wta)$$
\end{lem}
\begin{proof}
  We have
  \begin{align*}
    \label{}
    \qd(\wtb)\chmap_{\wtb+\rho}\big(\sch(V(\wta))\big) &=
    \chmap_{\rho}(L'_{0})^{-1} \chmap_{\wtb +\rho} \left(\sum_{w \in \weyl}
      \epsilon(w) \e^{w(\wta + \rho)}\right) \\
    & = \chmap_{\rho}(L'_{0})^{-1} \sum_{w \in \weyl}\epsilon(w)q^{2<w(\wta
      +\rho), \wtb +\rho>}\\
    & = \chmap_{\rho}(L'_{0})^{-1} \sum_{w \in \weyl}\epsilon(w)q^{2<w(\wtb
      +\rho), \wta +\rho>}\\
    & = \qd(\wta)\chmap_{\wta+\rho}\big(\sch(V(\wtb))\big)
  \end{align*}
  where the third equality follows from the fact that the form $<.,.>$ is
  symmetric and $\weyl$-invariant on $\h^{*}$.
\end{proof}

\begin{rk}
  One can easily check that Proposition \ref{P:S'} is also true if one
  replaces $V(\wta)$ or $V(\wtb)$ with its odd analog $V(\wta)^-$ or
  $V(\wtb)^-$, respectively.  As $\sch(V(\wta)^-)=-\sch(V(\wta))$, the
  appropriate extension of $\qd$ for ``odd'' modules is
  $\qd(V(\wta)^-)=-\qd(V(\wta))$. Hence Lemma \ref{L:sym} is still valid for
  ``odd'' modules.
\end{rk}
%
%
The following lemma shows that there is a typical $\g$-module whose tensor
product with itself is multiplicity free.  This implies that the endomorphism
ring of this tensor product is commutative, which we use to prove Lemma
\ref{L:Hopf2}.

\begin{lem}\label{L:VnoMulti}
  There exists a typical weight module $V(\wtc)$ such that $V(\wtc)\otimes
  V(\wtc)$ splits as a direct sum of irreducible typical modules with no
  multiplicity.
\end{lem}
\begin{proof}
  Let $\alpha$ be an irrational number.  Let $\wtc$ be the weight determined
  by $\wtc(h_{s})=\alpha$ and $\wtc(h_{i})=0$ for $i \neq s$.  We will show
  that $V(\wtc)$ has the desired properties.

  Let $\chi_{1}=\prod_{\alpha \in \roots^{+}_{\p1}}(1+\e^{-\alpha})$.  From
  Proposition \ref{P:formSch} we have that the character of the typical
  $\g$-module $V(\wta)$ is of the form
\begin{equation}
  \label{E:sch3}
  \ch(V(\wta))= \chi_{1}\chi_{0}(\wta)
\end{equation}
where $\chi_{0}(\wta)$ is the character of an even finite-dimensional $\g_{\p
  0}$-module.  For $V(\wtc)$ we have $\chi_{0}(\wtc)=\e^{\wtc}$.  Therefore,
$\ch(V(\wtc)\otimes V(\wtc))= \chi_{1}(\chi_{1}\e^{2\wtc})$.  From equation
\eqref{E:sch3} it is enough to show that $\chi_{1}\e^{2\wtc}$ is a sum of
characters of distinct irreducible $\g_{\p0}$-modules. (The corresponding
$\g$-module will be typical because the weights that appear in
$\chi_{1}(\chi_{1}\e^{2\wtc})$ take irrational values on $h_s$.)

In the following two cases the roots are expressed in terms of the standard
orthogonal basis $(\epsilon_{i})$ of $\h^{*}$ (see the Appendix).

\smallskip
\noindent
\textit{Case 1:} $\g=\slmn$.
Set $\delta_{i}=\epsilon_{i+m}$.  Recall that on $\h$,
$\sum\epsilon_i=\sum\delta_j$.  The set of odd positive roots is given by
$\roots_{\p 1}^{+}=\{\epsilon_{i}-\delta_{j}: 1 \leq i \leq m, 1 \leq j \leq
n\}$.  Set $x_{i}=\e^{\epsilon_{i}}$ and $y_{j}=\e^{\delta_{j}}$, then
$$\chi_1=\prod_{\alpha \in \roots^{+}_{\p1}}(1+\e^{-\alpha})=
\prod_{i,j}(1+y_{j}/x_{i})=(x_{1}x_{2}\cdots
x_{m})^{-n}\prod_{i,j}(x_{i}+y_{j}).$$
Let $s_{\lambda}(x)=s_{\lambda}(x_{1},\ldots,x_{m})$ be the Schur function
associated to a partition $\lambda$.  From the Cauchy identity we have
\begin{equation}
  \label{E:Cauchy}
  (x_{1}x_{2}\cdots
  x_{m})^{-n}\prod_{i,j}(x_{i}+y_{j})=(x_{1}x_{2}\cdots
  x_{m})^{-n}\sum_{\lambda\subset (n^m)}s_{\lambda}(x)s_{\hat{\lambda}'}(y)
\end{equation}
where the sum is over all partitions
$\lambda=(\lambda_{1},\ldots,\lambda_{m})$ such that $\lambda_{1}\leq n$, the
complementary partition
$\hat{\lambda}=(\hat\lambda_{1},\ldots,\hat\lambda_{m})$ is defined by
$\hat\lambda_{i}= n - \lambda_{m+1-i}$, and $\hat{\lambda}' $ is the conjugate
of $\hat\lambda$.

The Lie algebra $\g_{\p 0}$ is isomorphic to
$\sll(m)\times\sll(n)\times{\mathbb T}_1$.  The central element acts by a
complex number on any irreducible $\g_{\p 0}$-module.  This correspond to the
$\CC$-grading of $\ZZ[\weights]$.  So the sum \eqref{E:Cauchy} splits as a sum
of characters of $mn+1$ distinct $\g_{\p 0}$-modules following the spectral
decomposition of the central element of $\g_{\p 0}$~:
$$\e^{2\wtc}\chi_1=\sum_{i=0}^{mn}\left(\e^{2\wtc}(x_{1}x_{2}\cdots
x_{m})^{-n}
\sum_{\mbox{\scriptsize$\begin{array}{c}|\lambda|=i\\\lambda\subset
(n^m)\end{array}$}} s_{\lambda}(x)s_{\hat{\lambda}'}(y)\right)$$ 
Now for fixed $i$, the functions $s_{\lambda}(x)$ with
$|\lambda|=\lambda_1+\cdots+\lambda_m=i$ are the characters of distinct
irreducible $\sll(m)$-modules (with highest weight
$(\lambda_1-\lambda_2,\ldots,\lambda_{m-1}-\lambda_m)$).  Furthermore
$s_{\hat{\lambda}'}(y)$ is also the character of an irreducible
$\sll(n)$-module so all the terms in this sum are characters of different
irreducible $\g_{\p 0}$-modules.

\smallskip
\noindent
\textit{Case 2:} $\g=\osp$.
Set $\delta_{i}=\epsilon_{i+2}$ for $i=1,\ldots,n$.  The set of odd positive
roots is given by $\roots_{\p 1}^{+}=\{\epsilon_{1}\pm \delta_{i}: 1 \leq i
\leq n\}$.  Set $x_{i}=\e^{\delta_{i}}$ and $y=\e^{-\epsilon_{1}}$.  Let
$e_{k}(z_{1},\ldots,z_{p})=e_{k}(z)$ be the $k$\textsuperscript{th} elementary
symmetric function in $p$ variables $z_{1},\ldots,z_{p}$.  We have
$$\chi_1=\prod_{i=1}^{n}(1+yx_{i})(1+y/x_{i})=
\sum_{k=0}^{2n}e_{k}(x_{1},\ldots,x_{n},1/x_{1},\ldots,1/x_{n})y^{k}.$$
The proof follows from the fact that $e_{2n-k}(x,x^{-1})=e_k(x,x^{-1})$,
$e_0=1$ and for $k=1\cdots n$, $e_{k}(x,x^{-1})=\sum_{i=0}^{\lfloor
  n/2\rfloor} \Gamma_{k-2i}$ where $\Gamma_{j}$ is the character of the
irreducible $\spn$-module with highest weight the $j$\textsuperscript{th}
fundamental weight. (see \cite{FH} section 24.2 page 407).
\end{proof}
\begin{lem}\label{L:Hopf2}
  Let $V(\wtcc)$ be the module of Lemma \ref{L:VnoMulti}.  Let $\vz$ be the
  $\Uh$-module which is the deformation of $V(\wtcc)$.  Then we have
  $$\left<\, \put(0,15){\scriptsize\hbox{$\wtcc$}}
    \hspace{1ex}\epsh{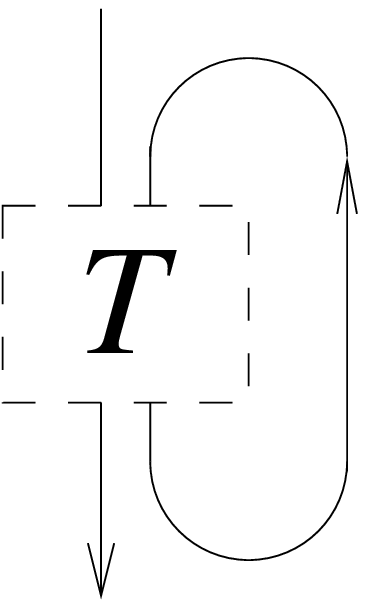}{8ex} \put(1,0){\scriptsize\hbox{$\wtcc$}}
    \hspace{1ex}\right> =\left< \put(-4,0){\scriptsize\hbox{$\wtcc$}}
    \hspace{1ex}\epsh{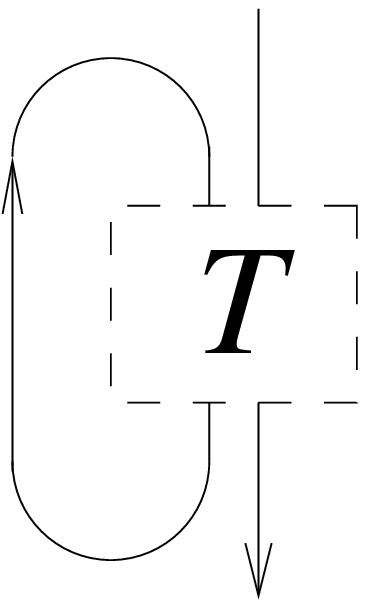}{8ex} \put(-4,15){\scriptsize\hbox{$\wtcc$}}
    \hspace{1ex}\right>$$ for all $T\in\T\big((\vz, \vz), (\vz, \vz)\big)$.
 \end{lem}
\begin{proof}
  Set $E=\End(\vz \otimes \vz)$.\\
  Consider the following linear forms on $E$:
  $$\tr_{L}(f)=(d_{\vz}\otimes \Id_{\vz})\circ(\Id_{\vz^{*}}\otimes
  f)\circ(b'_{\vz}\otimes \Id_{\vz}) \in \End(\vz)\cong \CC[[h]],$$
  $$\tr_{R}(f)=(\Id_{\vz}\otimes d'_{\vz}) \circ (f \otimes \Id_{\vz^{*}})
  \circ(\Id_{\vz}\otimes b_{\vz}) \in \End(\vz)\cong \CC[[h]].$$ 
  Let $T$ be any element of $\T((\vz, \vz),(\vz,\vz))$.  Let $F$ is the
  functor described in subsection \ref{SS:quant}.  In general,
  $$(\tr_{R}\circ F)(T)=(\tr_{L}\circ F)(c_{\vz,\vz}^{-1}Tc_{\vz,\vz})$$
  but from Lemma \ref{L:VnoMulti} we have that $E$ is a commutative algebra
  and so $F(c_{\vz,\vz}^{-1}Tc_{\vz,\vz})=F(T)$.
\end{proof}
%
%
%
The following lemma is a key ingredient in the proof of Theorem \ref{thF'}.
\begin{lem}\label{key}
  Let $\wtaa$ and $\wtbb$ be two typical weights.  Then we have
  $$\qd(\wtaa)\left<\, \put(4,15){\small\hbox{$\wtaa$}}
    \hspace{1ex}\epsh{fig37}{8ex} \put(2,0){\small\hbox{$\wtbb$}}
    \hspace{1ex}\right> =\qd(\wtbb)\left< \put(-2,0){\small\hbox{$\wtaa$}}
    \hspace{1ex}\epsh{fig38}{8ex} \put(-4,15){\small\hbox{$\wtbb$}}
    \hspace{1ex}\right>$$

  for all $T\in\T\Big(\big(\qum{V}(\wtaa),\qum{V}(\wtbb)\big),
  \big(\qum{V}(\wtaa),\qum{V}(\wtbb)\big)\Big)$.
\end{lem}
\begin{proof}
  Let $T\in\T\Big(\big(\qum{V}(\wtaa),\qum{V}(\wtbb)\big),
  \big(\qum{V}(\wtaa),\qum{V}(\wtbb)\big)\Big)$.  Let $V(\wtcc)$ be the module
  of Lemma \ref{L:VnoMulti}.  By definition we have
  \begin{align}
    \label{E:pict1}
    \left<\put(-1,15){\small\hbox{$\wtcc$}} \put(15,-14){\small\hbox{$\wtaa$}}
      \hspace{2ex}\epsh{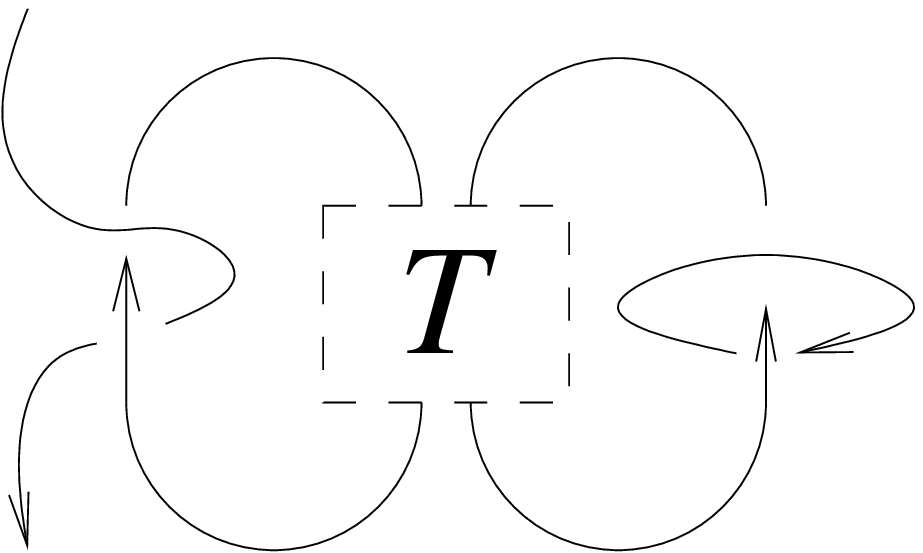}{8ex} \put(-10,-14){\small\hbox{$\wtbb$}}
      \put(2,2){\small\hbox{$\wtcc$}} \hspace{2ex}\right>&=
    \left< \put(4,15){\small\hbox{$\wtcc$}} \hspace{1ex}\epsh{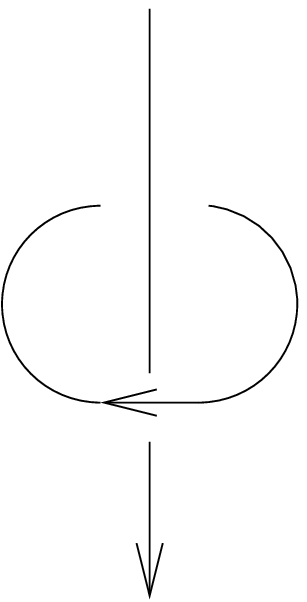}{8ex}
      \put(2,0){\small\hbox{$\wtaa$}} \hspace{1ex}\right>
    \left< \put(4,15){\small\hbox{$\wtaa$}} \hspace{1ex}\epsh{fig37}{8ex}
      \put(2,0){\small\hbox{$\wtbb$}} \hspace{1ex}\right>
    \left< \epsh{fig43}{8ex}\put(2,0){\small\hbox{$\wtcc$}}
      \put(-8,15){\small\hbox{$\wtbb$}}
      \hspace{2ex}\right> \notag\\
    &=S'(\wtaa,\wtcc)S'(\wtcc,\wtbb)\left<\, \put(4,15){\small\hbox{$\wtaa$}}
      \hspace{1ex}\epsh{fig37}{8ex} \put(2,0){\small\hbox{$\wtbb$}}
      \hspace{1ex}\right>.
  \end{align}

  Similarly,
  \begin{align}
    \label{E:pict2}
    \left<\put(-1,2){\small\hbox{$\wtcc$}} \put(15,-14){\small\hbox{$\wtaa$}}
      \hspace{2ex}\epsh{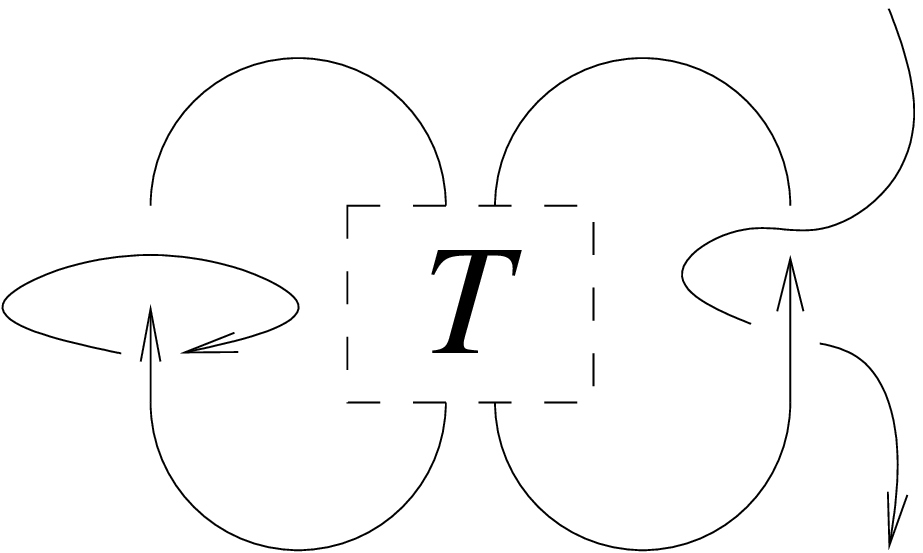}{8ex} \put(-10,-14){\small\hbox{$\wtbb$}}
      \put(1,15){\small\hbox{$\wtcc$}} \hspace{2ex}\right> &=
    S'(\wtbb,\wtcc)S'(\wtcc,\wtaa)\left< \put(-2,0){\small\hbox{$\wtaa$}}
      \hspace{1ex}\epsh{fig38}{8ex} \put(-4,15){\small\hbox{$\wtbb$}}
      \hspace{1ex}\right>.
  \end{align}

  From Lemma \ref{L:Hopf2} we have that the left sides of equations
  \eqref{E:pict1} and \eqref{E:pict2} are equal.  Thus, the results follows
  from Lemma \ref{L:sym} and Corollary \ref{C:SnotZero}.
\end{proof}
\begin{rk} \label{qdim} Let $\qum V_\wta$ be the deformation of a typical
  $\g$-module of highest weight $\wta$.  Compare the definition of $\qd$ given
  in \ref{D:d} with the quantum dimension of the $\Uhp$-module $V_0$ with the
  same highest weight $\wta$:
  $$\qdim(V_0)=\frac {\chmap_{\wta+\rho_{\p0}}(L'_{0})}
  {\chmap_{\rho_{\p0}}(L'_{0})} =\frac {\chmap_{\wta+\rho}(L'_{0})}
  {\chmap_{\rho}(L'_{0})}\quad\text{ so }\quad \qd(\qum V_\wta)= \frac
  {\qdim(V_0)} {\chmap_{\wta+\rho}(L'_{1})}.$$
  Also, compare $\qd(\qum V_\wta)$ with the formula that gives the quantum
  dimension of $\qum V_\wta$:
  $$\qdim(\qum V_\wta)=\chmap_{\rho}(\sch(\qum V_\wta)) =
  \chmap_{\rho}(L'_1)\frac {\chmap_{\wta+\rho}(L'_{0})}
  {\chmap_{\rho}(L'_{0})}$$ this last product vanishes because
  $\chmap_{\rho}(L'_1)=0$

  Lemma \ref{key} suggests to consider $\qd$ as a replacement for the quantum
  dimension.  Also, remark that any map proportional to $\qd$ would also be
  appropriate.  For example $h^{k_{\g}}\qd$ which admit a classical limit
  ($h\rightarrow0$).  Hence one gets a classical analog of Lemma \ref{key}
  ($\lim_{h\rightarrow0}h^{k_{\g}}\qd$ is a ``replacement'' for the
  super-dimension) which seems unknown.  In fact we can show, using the
  Kontsevich integral, that this classical version would be equivalent to
  Lemma \ref{key}.  But surprisingly, we have not found a simpler
  ``classical'' proof of this result.
\end{rk}

\begin{proof}[Proof of Theorem  \ref{thF'}]
  Any closed ribbon graph $L\in\T(\emptyset,\emptyset)$ over $\M$ with at
  least one edge colored by a typical module $\qum{V}(\wtaa)$ can be
  represented as the closure of $T_\wtaa\in\T(\qum{V}(\wtaa),\qum{V}(\wtaa))$.
  We set
  $$F'(L)=\qd(\wtaa)<T_\wtaa>$$.  If $L$ can also be represented as the
  closure of $T_\wtbb\in\T(\qum{V}(\wtbb),\qum{V}(\wtbb))$ for some typical
  weight $\wtbb$ then there exits
  $T\in\T\Big(\big(\qum{V}(\wtaa),\qum{V}(\wtbb)\big),
  \big(\qum{V}(\wtaa),\qum{V}(\wtbb)\big)\Big)$ such that $T_\wtaa=\left<\,
    \put(5,13){\scriptsize\hbox{$\wtaa$}} \hspace{1ex}\epsh{fig37}{6ex}
    \put(2,0){\scriptsize\hbox{$\wtbb$}} \hspace{1ex}\right>$ and
  $T_\wtbb=\left< \put(-1,0){\scriptsize\hbox{$\wtaa$}}
    \hspace{1ex}\epsh{fig38}{6ex} \put(-3,13){\scriptsize\hbox{$\wtbb$}}
    \hspace{1ex}\right>$ so by Lemma \ref{key} the definition of $F'(L)$ does
  not depend on the choice of $T_\wtaa$.
\end{proof}

\section{The multivariable invariant}\label{S:MVA}

In this section we construct families of multivariable invariants of links and
prove Theorem \ref{T:Mpoly2}.

\begin{lem}\label{L:twist}
  Let $\wta$ be a dominant weight.  Then the value of the twist
  $\twist_{\qum{V}(\wta)}$ is $q^{<\wta,\wta+2\rho>}$.  In other words,
  $$\left<
    \put(-2,0){\scriptsize\hbox{$\lambda$}}\epsh{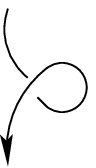}{6ex}\right> =
  q^{<\wta,\wta+2\rho>}.$$
\end{lem}
\begin{proof}
  The proof follow from Fact 1 in the proof of Proposition \ref{P:S'} and
  equation \eqref{E:K2}.
\end{proof}

Let $c,d\in \NN^{r-1}$.  Let $\wto a^c$ be the weight corresponding to an
$(r-1)$-tuple $c=(c_1,\ldots,c_{r-1})\in\NN^{r-1}$ and a complex parameter
$a\in\CC$.  Let $\V a^c$ be the $\Uh$-module $\qum{V}(\wto a^c)$.  Recall that
$\CC\setminus \TT_c$ is finite where $ \TT_c$ is the set of complex numbers
$a$ such that $\wto a^c$ is typical.

\begin{lem}\label{L:Mc}
  Recall the definition of $\qd$ given in Definition \ref{D:d}.  Let $\wto
  a^{c}$ be typical.  There exists integers $n^c_\alpha$ for
  $\alpha\in\roots^+$ such that for
  \begin{align*}
    M^c_0(q)&:=\prod_{\alpha \in
      \roots^{+}_{\p0}}\frac{q^{n^c_\alpha}-q^{-n^c_\alpha}}
    {q^{n^0_\alpha}-q^{-n^0_\alpha}} \in\ZZ[q^{\pm1}]\\
    \text{ and }\quad M^c_1(q,q_1)&:=\prod_{\alpha \in
      \roots^{+}_{\p1}}(q_1q^{n^c_\alpha}-q_1^{-1}q^{-n^c_\alpha})
    \in\ZZ[q^{\pm1},q_1^{\pm1}]
  \end{align*}
  one has $\chmap_{\wto a^{c}+\rho}(L'_{0})/\chmap_{\rho}(L'_{0})
  =M^c_0(\e^{h/2})$ and $\chmap_{\wto
    a^{c}+\rho}(L'_{1})=M^c_1(\e^{h/2},\e^{ha/2})$. In particular,
$$\qd(\wto a^{c})= \frac {M^c_0(\e^{h/2})}
{M^c_1(\e^{h/2},\e^{ah/2})}.$$
\end{lem}
\begin{proof}
  The existence of the integers $n^c_\alpha$ follows from an explicit
  computation (given in the Appendix) of the products $<\wto
  a^c+\rho,\alpha>$.  As $M^c_0(q)$ (see Remark \ref{qdim}) is the quantum
  dimension of the $\Uhp$-module with the highest weight $\wto a^{c}$, it is
  in $\ZZ[q^{\pm1}]$.
\end{proof}

In \cite{Le} Le proved that the quantum invariants arising from simple Lie
algebras (with a suitable normalization) are Laurent polynomial in $q$.  His
proof uses Lusztig's canonical basis.  We will now show that a similar
statement holds for $\g$.  However in our case, the complex parameter of
$\V{a}^{c}$ will lead to a Laurent polynomial in more than one variable.

\begin{lem}\label{L:R-mat}
  There exists bases $B_a^c$ of $\V{a}^c$ and elements
  $R^{c,d}(x,y,z)\in\GL(\ZZ[x^{\pm1},y^{\pm1},z^{\pm1}])$ such that the action
  of the $R$-matrix on $\V{a}^c\otimes \V b^d$ in the basis $B_a^c\times
  B_b^d$ is given by $q^{<\wto a^c, \wto b^d>}R^{c,d}(q,q^a,q^b)$.
\end{lem}
\begin{proof}
  We fix $a\in\CC\setminus \QQ$ and consider the ring of Laurent polynomials
  in two variables $A=\ZZ[q^{\pm1},q^{\pm a}]\subset\CC[[h]]$.  Let
  $\Uz^{\p0}$ (resp. $\Uz$) be the $A$-sub-algebra of $\Uhp$ (resp. of $\Uh$)
  generated by $X:=\{E_i^{(k)},F_i^{(k)},K_i,K_s : k\in\NN, i=1\cdots r, i\neq
  s\}$ (resp. by $X\cup\{F_s,(q-q^{-1})E_s\}$) where $G^{(k)}:=G^k/(k)_q!$ and
  $K_{i}:=q^{h_{i}}$.

  Recall that the $R$-matrix of $\Uh$ is of the form $R=\check{R}K$ where
  $\check{R}=\prod_{\alpha \in\roots^{+}}\check{R}_{\alpha}$.  It follows from
  Lusztig (\cite{L}, see also \cite{J}) that the quasi R-matrix of $\Uhp$,
  $\check R_{\p0}:=\prod_{\alpha \in \roots^{+}_{\p0}} R_{\alpha}$
  is an element of (a completion of) $\Uz^{\p0}\otimes_A\Uz^{\p0}$.  By
  definition $R_{\alpha}=1 + (q-q^{-1})E_{\alpha}\otimes F_{\alpha}$ for
  $\alpha\in \roots^{+}_{\p1}$.  Combining the last two sentences we have that
  $\check{R}$ is an element of (a completion of) $\Uz\otimes_A\Uz$.

  Let $v$ be a highest weight vector of $\V{a}^c$ and set $L=\Uz v$.  We will
  show that there exists a basis $B_a^c\times B_b^d$ such that the entries of
  the endomorphism of $\V{a}^c\otimes \V{b}^d$ arising from $\check{R}$ are in
  $\ZZ[q^{\pm1},q^{\pm a}, q^{\pm b}]$.  From the previous paragraph it
  suffices to show that there exists a basis for $L$ such that the action of
  any element of $\Uz$ on $L$ has a matrix (in this basis) whose entries are
  Laurent polynomials in $q$ and $q^{a}$.

Let $B_{\p0}$ be a Lusztig's canonical basis of
$L_{\p0}=\Uz^{\p0}v=\bigoplus_{x\in B_{\p0}}A.x$.  Let $B_{\p1}$ be the PBW
basis of $U_h\nil_{\p1}^-$ given by the $2^{|\roots_{\p1}^+|}$ ordered product
of elements from any subset of $\{F_\alpha:\alpha\in\roots_{\p1}^+\}$ (see
\cite{ZA,ZC,Dw}).  Let $B=\{f.x:f\in B_{\p1},\,x\in B_{\p0}\}$.  It follows
from the PBW theorem that $B$ generates $\V{a}^c$ (because $\Uh=
U_h\nil_{\p1}^- U_h\nil_{\p0}^-U_h\h U_h\nil_{\p0}^+ U_h\nil_{\p1}^+$ and
$\V{a}^c=\Uh v$).  Then the character formula for typical modules gives that
$B$ is a $\CC[[h]]$-basis for $\V{a}^c$.

We claim that $L=\bigoplus_{y\in B}A.y$ and hence $B_a^c:=B$ is a basis of
$L$.  This follows from the following two observations.  First, since $B$ is a
$\CC[[h]]$-basis for $\V{a}^c$ we have that $B_a^c$ is a linearly independent
subset of $L$.  Second, let $F=\bigoplus_{\alpha\in\roots^{+}_{\bar 1}}
A.F_\alpha\subset\Uz$ and $E=\bigoplus_{\alpha\in\roots^{+}_{\bar 1}}
A.E_\alpha\subset\Uz$.  The commutation relations of $\Uh$ (see
\cite{ZA,ZC,Dw}) imply that
\begin{align*}\label{E:com}
  \Uz^{\p0}.F&=F.\Uz^{\p0} & \Uz^{\p0}.E&=E.\Uz^{\p0}.
\end{align*}
From this equation and the commuting relation between $E_{\alpha}$ and
$F_{\beta}$, for $\alpha, \beta\in \roots_{\p 1}^+$, we have that any $g\in
\Uz$ can be written as an element in $\Uz^{\p 1, -}\Uz^{\p 0}\Uz^{\p1, +}$
where $\Uz^{\p 1, \pm}$ is generated by $\{(q-q^{-1})E_{\alpha}: \alpha \in
\roots_{\p 1}^+ \}$ or $\{ F_{\alpha}: \alpha \in \roots_{\p 1}^+ \}$,
respectively.  In other words, for any $g\in \Uz$ the element $gv$ can be
written as an $A$-linear combination of elements in $B_{a}^{c}$ any $g\in
\Uz$.  Hence, the free $A$-lattice generated by $B$ is stable by any $u\in\Uz$
and $\check R$ acts as desired.
Notice that all of the above could have be done over $\ZZ[q^{\pm1}]$
independently of the complex parameter $a$ except that some of the commuting
relations involve Laurent polynomials in $K_s$.  The element $K_s$ acts (by
multiplication) on the weight vectors of $\V{a}^c$ by $q^aq^k$ (for some
$k\in\ZZ$).  So the action of any element of $\Uz$ has a matrix in the basis
$B_a^c$ whose coefficients are Laurent polynomials independent of $a$ in the
two variables $q$ and $q^a$.

Finally we show that $K$ acts appropriately.  A vector in the basis $B_a^c$ is
a weight vector.  Its weight differs from the highest weight by an element of
the root lattice.  But for a root $\alpha$, one has $<\alpha,\wto
a^c>\in\ZZ+a\ZZ$.  This with equation \eqref{E:K1} shows that $q^{-<\wto a^c,
  \wto b^d>}K$ acts by monomials of $\ZZ[q^{\pm1},q^{\pm a},q^{\pm b}]$ on the
elements of the bases $B_{a}^c\times B_b^d$.
\end{proof}

\begin{rk}
  One has $<\wto a^c,\wto b^d>=\frac{n_1ab+n_2a+n_3b+n_4}{n_5}$ where $n_i\in
  \NN$. Thus, without the linking correction of Theorem \ref{T:MCpoly2}, $F'$
  takes values in a more complicated ring than a ring of Laurent polynomials.
\end{rk}

\begin{theo}\label{T:MCpoly2}
  Let $L$ be a link with $k$ ordered components colored by $k$ $(r-1)$-tuples
  $\bc=(\bc_1,\ldots,\bc_k)\in\left(\NN^{r-1}\right)^k$. Then there exists a
  multivariable Laurent polynomial $M(L;\bc)$ with values in
  $\ZZ[q^{\pm1},q_1^{\pm1},\ldots,q_k^{\pm1}]$ if $k\geq2$ and in
  $(M^{\bc_1}_1(q,q_1))^{-1} \ZZ[q^{\pm1},q_1^{\pm1}]$ if $k=1$ such that if
  \begin{itemize}
  \item $L'$ is any framed representative of $L$,
  \item $(\cn_1,\ldots,\cn_k)\in\TT_{\bc_1}\times\cdots\times\TT_{\bc_k}$,
  \item the $i$\textsuperscript{th} component of $L'$ is colored by the
    typical module $\qum{V}_{\cn_i}^{\bc_{i}}$ with highest weight
    $\wto{\cn_i}^{\bc_{i}}$,
  \end{itemize}
  then one has
  $$F'(L')=\e^{\sum lk_{ij}<\wto{\cn_i}^{\bc_{i}},
    \wto{\cn_j}^{\bc_{j}}+2\rho>h/2} {M(L;\bc)|_{q_i=\e^{\cn_ih/2}}}$$ where
  $lk_{ij}$ are the linking numbers of $L'$ (when $i=j$ the $lk_{ii}$ is the
  framing of $L'_i$).
\end{theo}
\begin{proof}
  Choose $k$ complex numbers $\cn_1,\ldots,\cn_k$ such that
  $(1,\cn_1,\ldots,\cn_k)$ is a linearly independent family of the
  $\QQ$-vector space $\CC$.  Let $\phi:
  \ZZ[q^{\pm1},q_1^{\pm1},\ldots,q_k^{\pm1}] \rightarrow\CC[[h]]$ be the ring
  map defined by
  $$\phi(q)=\e^{\frac{h}{2}} \; \text{ and }  \;
  \phi(q_i)=\e^{\frac{\cn_ih}{2}}$$ 
  Then $\phi$ is injective since the family $\{\phi(q^{t_0}q_1^{t_1}\cdots
  q_n^{t_n}): (t_0,\ldots ,t_n)\in\ZZ^{n+1}\}$ is free.

  Consider a $(1,1)$-tangle $T$ obtained by opening the
  $t$\textsuperscript{th} component of $L'$ for $t\in\{1\cdots k\}$.  By
  definition $F'(L')=\qd(\wto{\cn_t}^{\bc_{t}})<T>$.  Combining this with
  Lemmas \ref{L:R-mat} and \ref{L:Mc} we have
  \begin{equation}
    \label{E:F'}
    F'(L')= \e^{\sum
      lk_{ij}<\wto{\cn_i}^{\bc_{i}},\wto{\cn_j}^{\bc_{j}}>h/2}
    \frac{\phi\big(M_0^{\bc_t}(q)\big)}{\phi\big(M_1^{\bc_t}(q,q_t)\big)}
    \operatorname{Im}(\phi).
  \end{equation}

  As $<\wto{\cn_i}^{\bc_i},2\rho>\in\ZZ+\cn_i\ZZ$, (see Appendix) we can
  define
  \begin{equation}
    \label{D:Mc}
    M(L;\bc):=\frac{\phi^{-1} \left(\e^{-\sum
          lk_{ij}<\wto{\cn_i}^{\bc_i},\wto{\cn_j}^{\bc_j}+2\rho>h/2}
        M_1^{\bc_1}(\e^{h/2},\e^{\cn_1h/2})F'(L')\right)}
    {M^{\bc_{1}}_1(q,q_1)}.
  \end{equation}
  Note that Theorem \ref{thF'} implies that $M(L;\bc)$ is well defined.  In
  other words, we are always able to cut the first component of $L'$.  We
  added the $<\wto{\cn_i}^{\bc_i},2\rho>$ in equation \eqref{D:Mc} to make $M$
  a link invariant, i.e. framing independent (see Lemma \ref{L:twist}).

  Next we will show that if $k\geq 2$ then
  \begin{equation}
    \label{E:FinIm}
    F'(L')\in \e^{\sum
      lk_{ij}<\wto{\cn_i}^{\bc_{i}},\wto{\cn_j}^{\bc_{j}}>h/2}
    \operatorname{Im}(\phi).
  \end{equation}  
  For $i=1,2$, let $T_{i}$ be a $(1,1)$-tangles whose closure is $L'$ and
  whose open strand is the $i$th component of $L'$.  From Theorem \ref{thF'}
  we have $F'(T_{1})=F'(T_{2})$.  Then Lemma \ref{L:R-mat} implies the
  existence of Laurent polynomials $P_{1}$ and $P_{2}$ such that
  \begin{equation}
    \label{E:FPoly}
    \e^{-\sum lk_{ij}<\wto{\cn_i}^{\bc_{i}},\wto{\cn_j}^{\bc_{j}}>h/2}
    F'(L)= \qd(\wto{\cn_1}^{\bc_{1}})\phi(P_{1})=
    \qd(\wto{\cn_2}^{\bc_{2}}) \phi(P_{2}).
  \end{equation}
  By definition
  \begin{equation}
    \label{E:d}
    \qd(\wto{\cn_i}^{\bc_{i}})=\frac{M^{\bc_i}_0(\e^{h/2})}
    {M^{\bc_i}_1(\e^{h/2},\e^{h\wto{\cn_i}^{\bc_{i}}/2})},\;\; \text{ for
    } i=1,2.
  \end{equation}
  Thus
  $$M^{\bc_1}_0(q)M^{\bc_2}_1(q,q_{2})P_{1}=
  M^{\bc_2}_0(q)M^{\bc_1}_1(q,q_{1})P_{2}$$ 
  where $M^{\bc_1}_1(q,q_{1})$ and $M^{\bc_1}_0(q)M^{\bc_2}_1(q,q_{2})$ are
  relatively prime.  Since $\ZZ[q^{\pm1},q_1^{\pm1},\ldots,q_n^{\pm1}]$ is an
  unique factorization domain we have that $M^{c_1}_1(q,q_{1})$ divides
  $P_{1}$ and so equation \eqref{E:FinIm} holds.  Combining this with equation
  \eqref{D:Mc} we have
  $M(L;\bc)\in\ZZ[q^{\pm1},q_1^{\pm1},\ldots,q_n^{\pm1}]$.\\

  Because of Lemma~\ref{L:R-mat}, $M(L;\bc)$ is independent of the choice
  $\cn=(\cn_1,\ldots,\cn_k)$ lying in the dense subset of $\CC^n$ defined by
  the condition: $(1,\cn_1,\ldots,\cn_k)$ is a linearly independent family of
  the $\QQ$-vector space $\CC$.  Now the two maps $F'$ and $\phi\circ M$
  depend continuously of $\cn$ so the relation between $F'$ and $M$ in
  Theorem~\ref{T:MCpoly2} is valid for any
  $(\cn_1,\ldots,\cn_k)\in\TT_{\bc_1}\times\cdots\times\TT_{\bc_k}$.
\end{proof}
\begin{proof}[Proof of Theorem \ref{T:Mpoly2}]
  Just apply Theorem \ref{T:MCpoly2} to the link $L$ with all components
  colored by $c$.
\end{proof}
\section*{Appendix}

\subsection*{The $\sll(m|n)$ case}
%
$\sll(m|n)$ is the Lie superalgebra with Cartan matrix $A=\left(a_{ij}\right)$
whose non zeros entries given by
\begin{align*}
  a_{i,i}&=2&\text{ except }a_{m,m}=0\\
  a_{i,i+1}&=-1&\text{ except }a_{m,m+1}=1\\
  a_{i+1,i}&=-1.
\end{align*}
Set $d_i=1$ for $i=1\cdots m$ and $d_i=-1$ for $i>m$ then $(d_{i}a_{ij})$ is a
symmetric matrix.

We can identify $\sll(m|n)$ with the Lie superalgebra of super-trace zero
$(m|n)\times(m|n)$ matrices.  This standard representation is obtained by
sending $e_i$ on the elementary matrix $E_{i,i+1}$, $f_i$ on $E_{i+1,i}$,
$h_i$ on $E_{i,i}-E_{i+1,i+1}$ if $i\neq m$ and $h_m=E_{m,m}+E_{m+1,m+1}$. The
Cartan subalgebra $\h$ with basis $(h_i)$ is contained in the space of
diagonal matrices $X$.  The space $X^{*}$ has a canonical basis
$(\epsilon_1,\ldots\epsilon_{m+n})$ which is dual to the basis formed by the
matrices $E_{i,i}$. Set $\delta_{i}=\epsilon_{i+m}$, then $\h$ is the kernel
of the super-trace $\str=\sum\epsilon_i-\sum\delta_j$.  Therefore, $\h^*$ is
the quotient of $X^{*}$ by the super-trace.

The bilinear form on $\h$ given by $<h_i,h_j>=d_j^{-1}a_{ij}$ is equal to the
restriction on $\h$ of $<H,H'>=\str(H.H')$. So it extends to the whole set of
diagonal matrices and induces on its dual the bilinear form defined by
$<\epsilon_i,\epsilon_j>=\delta^i_j$ for $i,j=1\cdots m$ and
$<\delta_i,\delta_j>=-\delta^i_j$ for $i,j=1\cdots n$. Hence $\h^*$ can also
be identified as an euclidien space with $\str^\perp$.

The set of positive roots is $\roots^+=\roots^+_{\p 0}\cup\roots^+_{\p 1}$
with
$$\roots^+_{\p 0}=\{\epsilon_i-\epsilon_j,\,1\leq i<j\leq m\}\cup
\{\delta_i-\delta_j,\,1\leq i<j\leq n\}$$
$$\quad \text{and}\quad \roots^+_{\p 1}=\{\epsilon_i-\delta_j\}$$
The half sums of positive roots are given by
$$2\rho_{0}=\sum_i (m+1-2i)\epsilon_i+\sum_j (n+1-2j)\delta_j\quad
\text{,}\quad 2\rho_{1}=n\sum_i \epsilon_i-m\sum_j \delta_j$$
$$\text{and } \rho=\rho_{0}-\rho_{1}=\frac12\left(\sum_i
  (m-n+1-2i)\epsilon_i+\sum_j (m+n+1-2j)\delta_j\right).$$ 
Up to the super-trace, representatives of the fundamental weights are given by
\begin{align*}
  w_k&=\sum_{i=1}^k\epsilon_i\text{ for }k=1\cdots m-1\\
  w_m&=\sum_i \epsilon_i=\sum_j \delta_j\\
  w_{m+k}&=-\sum_{j=k+1}^{n}\delta_j\text{ for }k=1\cdots n-1.\\
\end{align*}
Thus for the weight $\wto a^c$ (with $c\in\NN^{r-1}$ and $a\in\TT_c$) of
section \ref{S:MVA}:
\begin{equation}\label{E:weight_sl}
  \wto a^c=c_1
  w_1+\cdots+c_{m-1}w_{m-1}+aw_m+c_mw_{m+1}+\cdots+c_{m+n-2}w_{m+n-1}.
\end{equation}
One has
\begin{align*}
  <\wto a^c,\rho_0>&=\frac1{2}\left(\sum_{i=1}^{m-1} i(m-i)c_i
    -\sum_{i=1}^{n-1}i(n-i)c_{m+n-1-i}\right)\\
  <\wto a^c,\rho_1>&=\frac1{2}\left(\sum_{i=1}^{m-1} nic_i +mna
    -\sum_{i=1}^{n-1}mic_{m+n-1-i}\right)\\
  <\wto a^c,2\rho>&=\sum_{i=1}^{m-1} i(m-n-i)c_i-mna
  +\sum_{i=1}^{n-1}i(m-n+i)c_{m+n-1-i}
\end{align*}
\begin{align*}
  <\epsilon_i-\delta_j,\wto a^c+\rho>&= <\epsilon_i-\delta_j,\wto
  a^c>&+& <\epsilon_i-\delta_j,\rho>\\
  &=\sum_{k=i}^{m-1}c_k+a-\sum_{k=m}^{m+n-1-j}c_k&+&(m+1-i-j)\\
  \text{and for }i<j,\\
  <\epsilon_i-\epsilon_j,\wto a^c+\rho>&= j-i+\sum_{k=i}^{j-1}c_k\\
  <\delta_i-\delta_j,\wto a^c+\rho>&= i-j-\sum_{k=m+n-i}^{m+n-j-1}c_{k}.
\end{align*}
(Remark that in all these scalar products, at least one of the two weights is
orthogonal to $\str$.)  This gives the form of $\qd(\wto a^c)$:
\begin{align*}
  \qd(\wto a^c)&=\frac{\chmap_{\wto a^c+\rho}(L'_{0})}{\chmap_{\wto a^c
      +\rho}(L'_{1})\chmap_{\rho}(L'_{0})}\\
  &=\prod_{\alpha\in\roots_{\p0}^+} \frac {q^{<\wto a^c+\rho,\alpha>}
    -q^{-<\wto a^c +\rho,\alpha>}} {q^{<\rho,\alpha>} -q^{-<\rho,\alpha>}}
  \Big/ \prod_{\alpha\in\roots_{\p1}^+}(q^{<\wto a^c +\rho,\alpha>}-q^{-<\wto
    a^c +\rho,\alpha>}).
\end{align*}

\subsection*{The $\osp$ case}
$\osp$ is the Lie superalgebra with Cartan matrix $A=\left(a_{ij}\right)$ with
non zeros entries given by
\begin{align*}
  a_{i,i}&=2&\text{ except }a_{1,1}=0\\
  a_{i,i+1}&=-1&\text{ except }a_{1,2}=1\text{ and }a_{n,n+1}=-2\\
  a_{i+1,i}&=-1.
\end{align*}
Set $d_1=1$, $d_i=-1$ for $i=2\cdots n$ and $d_{n+1}=-2$.  Then
$(d_{i}a_{ij})$ is a symmetric matrix.

Consider the super-symmetric form on $\CC^{2|2n}$ with matrix
$$B=\left(\begin{array}{ccc}
    I_2&0&0\\0&0&I_n\\0&-I_n&0\end{array}\right).$$
Then $\osp$ is the set of matrices of the form
$$X=\left(\begin{array}{rrr}A&B&C\\ -{^tC}&D&E\\ {^tB}&F&-{^tD}
  \end{array}\right)\text{ with }{^tA}=-A,\,{^tE}=E,\,{^tF}=F.$$
We consider the Cartan sub-algebra with basis $e_1=\sqrt{-1}(E_{1,2}-E_{2,1})$
and $e_{i+1}=E_{2+i,2+i}-E_{2+n+i,2+n+i}$ for $i=1\cdots n$. The dual basis is
denoted $(\epsilon,\delta_1,\ldots,\delta_n)$. The Cartan elements are
$h_i=e_i-e_{i+1}$ for $i=1\cdots n$ and $h_{n+1}=e_{n+1}$.

The bilinear form on $\h$ given by $<H,H'>=\frac12\str(H.H')$ induces a
bilinear form on $\h^*$ defined by $<\epsilon,\epsilon>=1$,
$<\epsilon,\delta_i>=0$ and $<\delta_i,\delta_j>=-\delta^i_j$ for $i,j=1\cdots
n$.

The set of positive roots is $\roots^+=\roots^+_{\p 0}\cup\roots^+_{\p 1}$
with
$$\roots^+_{\p 0}=\{\delta_i\pm\delta_j,\,1\leq i<j\leq n\}\cup
\{2\delta_i\}\quad \text{and}\quad \roots^+_{\p 1}=\{\epsilon\pm\delta_i\}$$
The half sums of positive roots are given by
$$\rho_{0}=\sum_i (n+1-i)\delta_i\quad
\text{,}\quad \rho_{1}=n\epsilon\quad \text{and}\quad
\rho=\rho_{0}-\rho_{1}.$$

The fundamental weights are given by
\begin{align*}
  w_1&=\epsilon\\
  w_{k+1}&=\epsilon+\sum_{i=1}^{k}\delta_i\text{ for }k=1\cdots n.
\end{align*}
Thus for the weight $\wto a^c$ (with $c\in\NN^{r-1}$ and $a\in\TT_c$) of
section \ref{S:MVA} :
\begin{equation}\label{E:weight_osp}
  \wto a^c=a w_1+c_1w_2+\cdots+c_nw_{n+1}.
\end{equation}
One has
$$<\wto a^c,\rho_1>=n(a+\sum_ic_i)\text{ and }
<\wto a^c,\rho_0>=-\sum_{i=1}^{n}\frac{i(2n-i+1)}2c_i$$
\begin{align*}
  <\epsilon\pm\delta_i,\wto a^c+\rho>&\begin{array}[t]{lcc}=
    <\epsilon\pm\delta_i,\wto a^c> &+& <\epsilon\pm\delta_i,\rho>\\\\
    =a+\sum_{k=1}^{n}c_k\mp\sum_{k=i}^{n}c_k &+&(-n\mp(n+1-i))
  \end{array}
  \\\text{ and for }i<j\\
  <\delta_i-\delta_j,\wto a^c+\rho>&= -\sum_{k=m+n-i}^{m+n-j-1}c_k+i-j.
\end{align*}
This gives the form of $\qd(\wto a^c)$:
\begin{align*}
  \qd(\wto a^c)&=\frac{\chmap_{\wto a^c+\rho}(L'_{0})}{\chmap_{\wto a^c
      +\rho}(L'_{1})\chmap_{\rho}(L'_{0})}\\
  &=\prod_{\alpha\in\roots_{\p0}^+} \frac {q^{<\wto a^c+\rho,\alpha>}
    -q^{-<\wto a^c +\rho,\alpha>}} {q^{<\rho,\alpha>} -q^{-<\rho,\alpha>}}
  \Big/ \prod_{\alpha\in\roots_{\p1}^+}(q^{<\wto a^c +\rho,\alpha>}-q^{-<\wto
    a^c +\rho,\alpha>}).
\end{align*}

\linespread{1}

\vfill


\begin{thebibliography}{99}
\setcounter{bibcount}{0}
\bibitem{ADO} Y. Akutsu, T. Deguchi, and T. Ohtsuki - {\em Invariants of
    colored links.} J. Knot Theory Ramifications \textbf{1} (1992), no. 2,
  161-184.

\bibitem{D} T. Deguchi - {\em Multivariable invariants of colored links
    generalizing the Alexander polynomial.} Proceedings of the Conference on
  Quantum Topology (Manhattan, KS, 1993) (River Edge, NJ), World Sci.
  Publishing, 1994, 67--85.

\bibitem{Dw} D. De Wit - {\em A Poincar\'e-Birkhoff-Witt commutator lemma for
    $U\sb q[{\rm gl}(m\vert n)]$.}  J. Math. Phys.  44 (2003), no. 1,
  315--327.

\bibitem{FH} W. Fulton, J. Harris - {\em Representation theory. A first
    course.} Graduate Texts in Mathematics, 129. Readings in Mathematics.
  Springer-Verlag, New York, 1991. xvi+551 pp.

\bibitem{G04B} N. Geer - {\em The Kontsevich integral and quantized Lie
    superalgebras.} Algebraic and Geometric Topology \textbf{5} (2005), paper
  no. 45, pages 1111-1139.
		
\bibitem{G05} N. Geer - {\em Some remarks on quantized Lie superalgebras of
    classical type.}  Preprint.

\bibitem{GP} N. Geer, B. Patureau-Mirand - {\em Multivariable link invariants
    arising from $\sll(2|1)$ and the Alexander polynomial.}  Preprint
  math.GT/0601291.

\bibitem{GP3} N. Geer, B. Patureau-Mirand - {\em On the Colored HOMFLY-PT,
    Multivariable and Kashaev Link Invariants.} In progress.

\bibitem{GPT} N. Geer, B. Patureau-Mirand, V. Turaev - {\em Renormalized
    quantum invariants.} In progress.

\bibitem{LG} M. Gould, J.R. Links - {\em Two variable link polynomials from
    quantum supergroups.} Lett. Math. Phys. \textbf{26} (1992), no. 3,
  187--198.

\bibitem{LGZ} M. Gould, J.R. Links, Y.Z. Zhang - {\em Type-I quantum
    superalgebras, $q$-supertrace, and two-variable link polynomials.}  J.
  Math. Phys.  \textbf{37} (1996), no. 2, 987--1003.

\bibitem{J} J.C. Jantzen - {\em Lectures on quantum groups.}  Graduate Studies
  in Mathematics, 6.  American Mathematical Society, Providence, RI, 1996.
  viii+266 pp.

\bibitem{K77} V.G. Kac - {\em Characters of typical representations of
    classical Lie superalgebras.} Comm. Algebra \textbf{5} (1977), no. 8,
  889--897.

\bibitem{K} V.G. Kac - {\em Lie superalgebras.} Advances Math. 26 (1977),
  8--96.

\bibitem{Kac78} V.G. Kac - {\em Representations of classical Lie
    superalgebras. } Differential geometrical methods in mathematical physics,
  II (Proc. Conf., Univ. Bonn, Bonn, 1977), pp. 597--626, Lecture Notes in
  Math., 676, Springer, Berlin, 1978.

\bibitem{Kv} R.M. Kashaev - {\em A link invariant from quantum dilogarithm.}
  Modern Phys. Lett. A \textbf{10} (1995), no. 19, 1409--1418.

\bibitem{KTol} S.M. Khoroshkin, V.N. Tolstoy - {\em Universal $R$-matrix for
    quantized (super)algebras.} {Comm. Math. Phys.}  \textbf{141} (1991), no.
  3, 599--617.

\bibitem{Le} T.T.Q. Le - {\em Integrality and symmetry of quantum link
    invariants.} Duke Math. J. 102 (2000), no. 2, 273--306.

\bibitem{LZ} J. R. Links, R. B. Zhang - {\em Multiparameter link invariants
    from quantum supergroups.}  J. Math. Phys. \textbf{35} (1994), no. 3,
  1377--1386.

\bibitem{L} G. Lusztig - {Introduction to quantum groups,} Birkhauser, Boston,
  1994.


\bibitem{Mur93} J. Murakami - {\em A state model for the multivariable
    Alexander polynomial.}  Pacific J. Math.  \textbf{157} (1993), no. 1,
  109--135.

\bibitem{MM} H. Murakami, J. Murakami - {\em The colored Jones polynomials and
    the simplicial volume of a knot.}  Acta Math. 186 (2001), no. 1, 85--104.

\bibitem{Tu} V.G. Turaev - {\em Quantum invariants of knots and 3-manifolds.}
  de Gruyter Studies in Mathematics, 18. Walter de Gruyter \& Co., Berlin,
  (1994).

\bibitem{Yam94} H. Yamane - {\em Quantized enveloping algebras associated with
    simple Lie superalgebras and their universal $R$-matrices.}  \emph{Publ.
    Res. Inst. Math. Sci.}  \textbf{30} (1994), no. 1, 15--87.

\bibitem{ZA} R.B. Zhang - {\em Finite-dimensional irreducible representations
    of the quantum supergroup $U_q({\rm gl}(m/n))$.} J. Math. Phys. 34 (1993),
  no. 3, 1236--1254.

\bibitem{ZC} R.B. Zhang - {\em Finite-dimensional representations of
    $U_q(C(n+1))$ at arbitrary $q$.} J. Phys. A 26 (1993), no. 23, 7041--7059.

\end{thebibliography}
\end{document}